\documentclass[11pt,reqno,a4letter]{article} 

\usepackage{amsthm,amsfonts,amscd,amsmath}
\usepackage[hidelinks]{hyperref} 
\usepackage{url}
\usepackage{pdfpages}
\usepackage[usenames,dvipsnames]{xcolor}
\hypersetup{
    colorlinks,
    linkcolor={red!70!white},
    citecolor={blue!70!white},
    urlcolor={blue!80!white}
}

\usepackage[normalem]{ulem}
\usepackage{graphicx} 
\usepackage{datetime} 
\usepackage{cancel}
\usepackage{verbatim}
\usepackage{subcaption}
\numberwithin{figure}{section}
\numberwithin{table}{section}

\usepackage{stmaryrd} 
\usepackage{framed} 
\usepackage{ulem}
\usepackage[T1]{fontenc}
\usepackage{pbsi}

\usepackage{enumitem}
\setlist[itemize]{leftmargin=0.65in}

\usepackage{rotating,adjustbox}

\usepackage{diagbox}

\let\citep\cite

\newcommand{\gkpSI}[2]{\ensuremath{\genfrac{\lbrack}{\rbrack}{0pt}{}{#1}{#2}}} 
\newcommand{\gkpSII}[2]{\ensuremath{\genfrac{\lbrace}{\rbrace}{0pt}{}{#1}{#2}}}
\newcommand{\cf}{\textit{cf.\ }} 
\newcommand{\Iverson}[1]{\ensuremath{\left[#1\right]_{\delta}}}

\newcommand{\seqnum}[1]{\href{http://oeis.org/#1}{\color{ProcessBlue}{\underline{#1}}}}

\usepackage{upgreek,dsfont}
\renewcommand{\chi}{\upchi}

\usepackage{ifthen}
\newcommand{\Hn}[2]{
     \ifthenelse{\equal{#2}{1}}{H_{#1}}{H_{#1}^{\left(#2\right)}}
}

\newcommand{\Floor}[2]{\ensuremath{\left\lfloor \frac{#1}{#2} \right\rfloor}}

\DeclareMathOperator{\DGF}{DGF}

\DeclareMathOperator{\mock}{mock}

\title{A catalog of interesting and useful Lambert series identities} 
\author{{\large Dr. Maxie Dion Schmidt} \\ 
        {\normalsize \href{mailto:maxieds@gmail.com}{maxieds@gmail.com}} \\ 
} 

\date{\small\underline{Last Revised:} \today\ \ -- \ \ Compiled with \LaTeX2e} 

\theoremstyle{plain} 
\newtheorem{theorem}{Theorem}

\numberwithin{theorem}{section}

\theoremstyle{definition}

\newtheorem{definition}[theorem]{Definition}
\newtheorem{notation}[theorem]{Notation}

\usepackage[top=15mm, bottom=15mm, left=12mm, right=12mm]{geometry}
\numberwithin{equation}{section}

\allowdisplaybreaks 

\begin{document} 

\maketitle

\begin{abstract} 
A \emph{Lambert series generating function} is a special series summed over an arithmetic 
function $f$ defined by 
\[
L_f(q) := \sum_{n \geq 1} \frac{f(n) q^n}{1-q^n} = \sum_{m \geq 1} (f \ast 1)(m) q^m. 
\]
Because of the way the left-hand-side terms of this type of generating function  
generate divisor sums of $f$ convolved by Dirichlet convolution with one, these expansions 
are natural ways to enumerate the ordinary generating functions of many multiplicative 
special functions in number theory. We present an overview of key properties of 
Lambert series generating function expansions, their more combinatorial generalizations, and 
include a compendia of tables illustrating known formulas for special cases of these series. 
In this sense, we focus more on the formal properties of the sequences that are enumerated by 
the Lambert series, and do not spend significant time treating these series as analytic objects 
subject to rigorous convergence constraints. 

The first question a reader 
might ask before examining this document is: 
\emph{Why was a catalog of interesting Lambert series identities compiled?} 
As with the indispensible reference by H. W. Gould and T. Shonhiwa, 
\emph{A catalog of interesting Dirichlet series}, 
for Dirichlet series (DGF) identities, there are many situations in which one needs a 
summary reference on Lambert series and their properties. 
New work has been authored recently tying Lambert series expansions to partition functions by expansions of 
their generating functions. 
In addition to these new expansions and providing an introduction to Lambert series, 
we have listings of classically relevant and ``odds and ends'' examples for 
Lambert series summations that are occasionally useful in applications. 

\bigskip
%\bigskip\hrule\bigskip
\textbf{Primary Math Subject Classifications (2010):} {\it 05A15; 11Y70; 11A25; and 11-00. }
\end{abstract}

\newpage
%\section{Notation} 
\label{Appendix_Glossary_NotationConvs}


# 1 Notation

| Symbol | Definition |
|---|---|
| $B_n(x)$, $B_n$ | The Bernoulli polynomials and Bernoulli numbers $B_n = B_n(0)$. These polynomials can be used via Faulhaber's formula, among others, to generate the integral $k^{th}$ power sums $\sum_{d|n} \phi_k(d)(n/d)^k = 1^k + 2^k + \cdots + n^k$. |
| $\binom{n}{k}$ | The binomial coefficients, $\binom{n}{k} = \frac{n!}{k!(n-k)!}$. |
| $\lceil x \rceil$ | The ceiling function $\lceil x \rceil := x + 1 - \{x\}$ where $0 \leq \{x\} < 1$ denotes the fractional part of $x \in \mathbb{R}$. |
| $[q^n]F(q)$ | The coefficient of $q^n$ in the power series expansion of $F(q)$ about zero. |
| $c_q(n)$ | Ramanujan's sum, $c_q(n) := \sum_{d|(q,n)} d\mu(q/d)$. |
| $D_z^{(j)}$ | The higher-order $j^{th}$ derivative operator with respect to $z$. |
| $d_k(n)$ | The generalized $k$-fold divisor function $d_k(n) = 1_{*_k}(n)$ whose DGF is $\zeta(s)^k$. Note that the divisor function $d(n) \equiv d_1(n)$. |
| $d(n)$ | The ordinary divisor function, $d(n) := \sum_{d|n} 1$. |
| $\varepsilon(n)$ | The multiplicative identity with respect to Dirichlet convolution, $\varepsilon(n) = \delta_{n,1}$. |
| $*; f * g$ | The Dirichlet convolution of $f$ and $g$, $f * g(n) := \sum_{d|n} f(d)g(n/d)$, for $n \geq 1$. This symbol for the discrete convolution of two arithmetic functions is the only notion of convolution of functions we employ within the article. |
| $f^{-1}(n)$ | The Dirichlet inverse of $f$ with respect to convolution defined recursively by $f^{-1}(n) = -\frac{1}{f(1)} \sum_{\substack{d|n \\ d>1}} f(d) f^{-1}(n/d)$ provided that $f(1) \neq 0$. |
| $\lfloor x \rfloor$ | The floor function $\lfloor x \rfloor := x - \{x\}$ where $0 \leq \{x\} < 1$ denotes the fractional part of $x \in \mathbb{R}$. |
| $f_{*_j}$ | Sequence of nested $j$-convolutions of an arithmetic function $f$ with itself for integers $j \geq 1$. We define $f_{*_0}(n) = \delta_{n,1}$, the multiplicative identity with respect to Dirichlet convolution. |
| $\gamma(n)$ | The squarefree kernel of $n$, $\gamma(n) := \prod_{p|n} p$. |
| $\gcd(m,n), (m,n)$ | The greatest common divisor of $m$ and $n$. Both notations for the GCD are used interchangably within the article. |
| $G_j$ | Denotes the interleaved (or generalized) sequence of pentagonal numbers defined explictly by the formula $G_j := \frac{1}{2} \left\lceil \frac{j}{2} \right\rceil \left\lceil \frac{3j+1}{2} \right\rceil$. The sequence begins as $\{G_j\}_{j \geq 0} = \{0, 1, 2, 5, 7, 12, 15, 22, 26, 35, 40, 51, \ldots\}$. |
| $\mathrm{Id}_k(n)$ | The power-scaled identity function, $\mathrm{Id}_k(n) := n^k$ for $n \geq 1$. The Dirichlet inverse of this function is given by $\mathrm{Id}_k^{-1}(n) := n^k \cdot \mu(n)$. |
| $\mathbb{1}_{\mathbb{S}}, \chi_{\text{cond}(x)}$ | We use the notation $\mathbb{1}, \chi : \mathbb{N} \to \{0, 1\}$ to denote indicator, or characteristic functions. In paticular, $\mathbb{1}_{\mathbb{S}}(n) = 1$ if and only if $n \in \mathbb{S}$, and $\chi_{\text{cond}}(n) = 1$ if and only if $n$ satisfies the condition cond. |
| $[n = k]_\delta$ | Synonym for $\delta_{n,k}$ which is one if and only if $n = k$, and zero otherwise. |
| $[\text{cond}]_\delta$ | For a boolean-valued cond, $[\text{cond}]_\delta$ evaluates to one precisely when cond is true, and zero otherwise. |
| $\vartheta_i(z, q), \vartheta_i(q)$ | For $i = 1, 2, 3, 4$, these are the classical Jacobi theta functions where $\vartheta_i(q) \equiv \vartheta_i(0, q)$. |

| Symbol | Definition |
|---|---|
| $J_t(n)$ | The Jordan totient function $J_t(n) = n^k \prod_{p|n}(1-p^{-t})$ also satisfies $\sum_{d|n} J_t(d) = n^t$. |
| $\Lambda(n)$ | The von Mangoldt lambda function $\Lambda(n) = \sum_{d|n} \log(d)\mu(n/d)$. |
| $\lambda(n)$ | The Liouville lambda function $\lambda(n) = (-1)^{\Omega(n)}$. |
| $\lambda_k(n)$ | The arithmetic function defined by $\lambda_k(n) = \sum_{d|n} d^k \lambda(d)$. |
| $L_f(q)$ | The lambert series generating function of an arithmetic function $f$, defined by $$L_f(q) := \sum_{n \geq 1} \frac{f(n)q^n}{1-q^n}, |q| < 1.$$ |
| $\text{lcm}(m,n), [m,n]$ | The least common multiple of $m$ and $n$. |
| $\log$ | The natural logarithm function, $\log(n) \equiv \ln(n)$. |
| $\text{lsb}(n)$ | The least significant bit of $n$ in the base-2 expansion of $n$. |
| $\mu(n)$ | The Möbius function, defined for $n \geq 1$ with $n = p_1^{\alpha_1} \cdots p_k^{\alpha_k}$ its factorization into distinct prime powers with $\alpha_i \geq 1$ for all $1 \leq i \leq k$ by $$\mu(n) = \begin{cases} 1, & \text{if } n=1; \\ (-1)^k, & \text{if } \alpha_i = 1, \forall 1 \leq i \leq k; \\ 0, & \text{otherwise.} \end{cases}$$ |
| OGF | *Ordinary generating function.* Given a sequence $\{f_n\}_{n \geq 0}$, its OGF (or sometimes called ordinary power series, OPS) enumerates the sequence by powers of a typically formal variable $z$: $F(z) := \sum_{n \geq 0} f_n z^n$. For $z \in \mathbb{C}$ within some radius or abcissa of convergence for the series, asymptotic properties can be extracted from the closed-form representation of $F$, and/or the original sequence terms can be recovered by performing an inverse $Z$-transform on the OGF. |
| $\omega_a$ | A primitive $a^{th}$ root of unity $\omega_a = \exp(2\pi \imath/a)$ for integers $a \geq 1$. |
| $\omega(n), \Omega(n)$ | If $n = p_1^{\alpha_1} \cdots p_r^{\alpha_r}$ is the prime factorization of $n$ into distinct prime powers, then $\omega(n) = r$ and $\Omega(n) = \alpha_1 + \cdots + \alpha_r$. |
| $\phi_k(n)$ | Generalized totient function, $\phi_k(n) := \sum_{\substack{1 \leq d \leq n \\ (d,n)=1}} d^k$. |
| $\phi(n)$ | Euler's classical totient function, $\phi(n) := \sum_{\substack{1 \leq d \leq n \\ (d,n)=1}} 1$. |
| $\Phi_n(z)$ | The $n^{th}$ cyclotomic polynomial in $z$ defined by $\Phi_n(z) := \prod_{\substack{1 \leq k \leq n \\ (k,n)=1}} (z - e^{2\pi \imath k/n})$. |
| $p(n)$ | The partition function generated by $p(n) = [q^n] \prod_{n \geq 1}(1-q^n)^{-1}$. |
| $\pi(x)$ | The prime counting function denotes the number of primes $p \leq x$. |
| $\sum_p, \sum_{p \leq x}$ | Unless otherwise specified by context, we use the index variable $p$ to denote that the summation is to be taken only over prime values within the summation bounds. |
| $\Psi_k(n)$ | The $k^{th}$ Dedekind totient function, $\Psi_k(n) := n^k \prod_{p|n}(1+p^{-k})$. |
| $\psi_k(n)$ | The arithmetic function defined by $\psi_k(n) := \sum_{d|n} d^k \mu^2(n/d)$. |

| Symbol | Definition |
|---|---|
| $(a;q)_n$, $(q)_n$ | The $q$-Pochhammer symbol defined as the product $(a;q)_\infty := \prod_{n \geq 1}(1-aq^{n-1})$. We adopt the notation that $(q)_n \equiv (q;q)_n$ and that $(a;q)_\infty$ denotes the limiting case for $|q| < 1$ as $n \to \infty$. |
| $(a_1,\ldots,a_r;q)_n$ | We use the common shorthand that $(a_1,\ldots,a_r;q)_n = \prod_{i=1}^{r}(a_i;q)_n$. |
| $r_k(n)$ | The sum of $k$ squares function denotes the number of integer solutions to $n = x_1^2 + \cdots + x_k^2$. A generating function is given by $r_k(n) = [q^n]\vartheta_3(q)^k$. |
| $\sigma_\alpha(n)$ | The generalized sum-of-divisors function, $\sigma_\alpha(n) := \sum_{d\mid n} d^\alpha$, for any $n \geq 1$ and $\alpha \in \mathbb{C}$. Note that the divisor function is also denoted by $d(n) \equiv \sigma_0(n)$. |
| $\left[{n \atop k}\right]$, $\left\{{n \atop k}\right\}$ | The Stirling numbers of the first and second kinds, respectively. Alternate notation for these triangles is given by $s(n,k) = (-1)^{n-k}\left[{n \atop k}\right]$ and $S(n,k) = \left\{{n \atop k}\right\}$. |
| $\tau(n)$ | The function defined by $\tau(n) := [x^{n-1}]\prod_{m\geq 1}(1-x^m)^{24}$. |
| $\zeta(s)$ | The Riemann zeta function, defined by $\zeta(s) := \sum_{n\geq 1} n^{-s}$ when $\Re(s) > 1$, and by analytic continuation to the entire complex plane with the exception of a simple pole at $s = 1$. |


\clearpage

%\verbatiminput{main.slo}
%\vskip -0.35in
%\printglossary[type={symbols},title={},style={glossstyleSymbol},nogroupskip=true]

\newpage
\section{Background and identities for Lambert series generating functions} 

\subsection{Introduction} 

In the most general setting, we define the \emph{generalized Lambert series} expansion 
for integers $0 \leq \beta < \alpha$ and any fixed arithmetic function $f$ as 
\[
L_f(\alpha, \beta; q) := \sum_{n \geq 1} \frac{f(n) q^{\alpha n - \beta}}{1-q^{\alpha n-\beta}}; 
     |q| < 1, 
\]
where the series coefficients of the Lambert series generating function are given by the divisor sums 
\[
[q^n] L_f(\alpha, \beta; q) = \sum_{\alpha d-\beta| n} f(d). 
\]
If we set $(\alpha, \beta) := (1, 0)$, the we recover the classical form of the Lambert series 
construction, which we will denote by the function 
$L_f(q) \equiv L_f(1, 0; q)$. 

There is a natural correspondence between a sequence's \emph{ordinary generating function} (OGF), and its 
Lambert series generating function. Namely, if $\widetilde{F}(q) := \sum_{m \geq 1} f(m) q^m$ is the OGF of 
$f$, then $$L_f(\alpha, \beta; q) = \sum_{n \geq 1} \widetilde{F}(q^{\alpha n-\beta}).$$ 
We also have a so-called \emph{Lambert transform} defined by \cite[\S 2]{RAMANUJAN-LSERIES-SURVEY}
\[
F_a(x) := \int_0^{\infty} \frac{t a(t)}{e^{xt}-1} dt. 
\]
This transform satisfies an inversion relation of the form 
\[
\tau a(\tau) = \lim_{k \rightarrow \infty} \frac{(-1)^{k+1}}{k!} \left(\frac{k}{\tau}\right)^{k+1} 
     \sum_{n \geq 1} \mu(n) n^k F_a^{(k)}\left(\frac{nk}{\tau}\right). 
\] 
The interpretation of this invertible transform as a so-called Lambert transformation considers taking 
the mappings $t a(t) \leftrightarrow a_n$ and $e^{-x} \leftrightarrow q$. 

\subsection{Higher-order derivatives of Lambert series} 

Ramanujan discovered the following remarkable identities \cite[\S 2]{RAMANUJAN-LSERIES-SURVEY}: 
\begin{subequations}
\begin{align}
\sum_{n \geq 1} \frac{(-1)^{n-1} q^n}{1-q^n} & = \sum_{n \geq 1} \frac{q^n}{1+q^n} \\ 
\sum_{n \geq 1} \frac{n q^n}{1-q^n} & = \sum_{n \geq 1} \frac{q^n}{(1-q^n)^2} \\ 
\sum_{n \geq 1} \frac{(-1)^{n-1} n q^n}{1-q^n} & = \sum_{n \geq 1} \frac{q^n}{(1+q^n)^2} \\ 
\sum_{n \geq 1} \frac{q^n}{n(1-q^n)} & = \sum_{n \geq 1} \frac{q^n}{1+q^n} \\ 
\sum_{n \geq 1} \frac{(-1)^{n-1} q^n}{n(1-q^n)} & = \sum_{n \geq 1} \log\left(\frac{1}{1-q^n}\right) \\ 
\sum_{n \geq 1} \frac{\alpha^n q^n}{1-q^n} & = \sum_{n \geq 1} \log\left(1+q^n\right) \\ 
\sum_{n \geq 1} \frac{n^2 q^n}{1-q^n} & = \sum_{n \geq 1} \frac{q^n}{(1-q^n)^2} \sum_{k=1}^{n} \frac{1}{1-q^k}. 
\end{align}
\end{subequations}
It follows that we can relate the partition function generating functions $(q; \pm q)_{\infty}$ to the 
exponential of the two logarithmically termed series above 
(\cf the remarks in Section \ref{subSection_OtherProps_RefsToq-SeriesExps}). 

More generally, higher-order $j^{th}$ derivatives for integers $j \geq 1$ can be obtained by differentiating 
the Lambert series expansions termwise in the forms of 
\citep[\cf Lemma 3]{MDS-COMBRESTRDIVSUMS-INTEGERS} 
\begin{subequations}
\begin{align}
q^j \cdot D_q^{(j)}\left[\frac{q^n}{1-q^n}\right] & = 
     \sum_{m=0}^j \sum_{k=0}^m \gkpSI{j}{m} 
     \gkpSII{m}{k} \frac{(-1)^{j-k} k! n^m}{(1-q^n)^{k+1}}, \\ 
     & = \sum_{r=0}^j \left[\sum_{m=0}^j \sum_{k=0}^m \gkpSI{j}{m}
     \gkpSII{m}{k} \binom{j-k}{r} \frac{(-1)^{j-k-r} k! n^m}{(1-q^n)^{k+1}}\right] q^{(r+1)n}. 
\end{align}
\end{subequations}
Here, since we can express the coefficients for all finite $n \geq 1$ as \[
[q^n] L_f(q) = [q^n] \sum_{m \leq n} \frac{f(m) q^m}{1-q^m}, 
\]
by partial sums of the generating functions, we need not worry about uniform convergence of the 
Lambert series generating function in $q$. 

By the binomial generating functions given by $[z^n] (1-z)^{m+1} = \binom{n+m}{m}$, we find that 
\[
[q^n] \left(\sum_{n \geq t} \frac{f(n) q^{mn}}{(1-q^n)^{k+1}}\right) = 
     \sum_{\substack{d|n \\ t \leq d \leq \Floor{n}{m}}} 
     \binom{\frac{n}{d} - m + k}{k} f(d), 
\]
for positive integers $m,t \geq 1$ and $k \geq 0$. 

\subsection{Relation of the coefficients of the classical Lambert series to Ramanujan sums} 

We define the functions $\widetilde{\Phi}_n(q)$ as the change of variable into the 
logarithmic derivatives of the \emph{cyclotomic polynomials} as 
\begin{align*} 
\widetilde{\Phi}_n(q) & = \frac{\Phi_n^{\prime}(1/q)}{q \cdot \Phi_n(1/q)} 
       = \frac{1}{q} \frac{d}{dw}\left[\sum_{d|n} \mu(n/d) \log(w^d-1)\right] \Biggr\rvert_{w=1/q} \\ 
     & = \sum_{d|n} \frac{d \mu(n/d)}{1-q^d}. 
\end{align*}
We can express the component series terms for 
$n \geq 1$ in the form of 
\citep{SCHMIDT-SODFORMULAS}
\[
\frac{1}{1-q^n} = \frac{1}{n} \sum_{d|n} \widetilde{\Phi}_d(q). 
\]
Then we can express the Lambert series coefficients, $(f \ast 1)(n)$, in terms of the \emph{Ramanujan sums}, $c_q(n)$, 
for each positive natural number $x \geq 1$ as 
\[
[q^x] \sum_{n \leq x} \frac{f(n)}{1-q^n} = 
     \sum_{d=1}^{x} c_d(x) \sum_{n=1}^{\Floor{x}{d}} \frac{f(nd)}{nd} = 
     \sum_{n \leq x} \frac{f(n)}{n} \sum_{d|n} c_d(x). 
\]

\subsection{Factorization theorems} 

\subsubsection{Classical series cases} 

The first form of the factorization theorems considered in \cite{AA,MERCA-LSFACTTHM} expands two variants of 
$L_f(q)$ as 
\[
\sum_{n \geq 1} \frac{f(n) q^n}{1 \pm q^n} = \frac{1}{(\mp q; q)_{\infty}} \sum_{n \geq 1} 
     \left(s_o(n, k) \pm s_e(n, k)\right) f(k) q^n,
\]
where $s_o(n, k) \pm s_e(n, k) = [q^n] (\mp q; q)_{\infty} \frac{q^k}{1 \pm q^k}$ is defined as the 
sum (difference) of the functions $s_o(n, k)$ and $s_e(n, k)$, which respectively denote the 
number of $k$'s in all partitions of $n$ into and odd (even) number of distinct parts. 
If we define $s_{n,k} = s_o(n, k) - s_e(n, k)$, then this sequence is lower triangular and invertible. 
Its inverse matrix is defined by \cite[\seqnum{A133732}]{OEIS} 
\[
s_{n,k}^{-1} = \sum_{d|n} p(d-k) \mu\left(\frac{n}{d}\right). 
\]
We can define the form of another factorization of $L_f(q)$ where $|C(q)| < \infty$ for all $|q| < 1$ is such that 
$C(0) \neq 0$ as 
\[
\sum_{n \geq 1} \frac{f(n) q^n}{1-q^n} = 
     \frac{1}{C(q)} \sum_{n \geq 1} \left(\sum_{k=1}^n s_{n,k}(\gamma) \widetilde{f}(k)(\gamma)\right) q^n,
\]
for any prescribed non-zero arithmetic function $\gamma(n)$ with 
\[
\widetilde{f}(k)(\gamma) = \sum_{d|k} \sum_{r| \frac{k}{d}} f(d) \gamma(r). 
\]
In this case, we have that 
\[
s_{n,k}^{-1}(\gamma) = \sum_{d|n} [q^{d-k}] \frac{1}{C(q)} \gamma\left(\frac{n}{d}\right). 
\]
This notion of factorization can be generalized to expanding the generalized Lambert series 
$L_f(\alpha, \beta; q)$ from the first subsection \cite{MERCA-SCHMIDT-RAMJ}. 

In either case, the coefficients generated by $L_f(q)$ as $[q^n] L_f(q) = (f \ast 1)(n)$ and their 
summatory functions, 
\[
\Sigma_f(x) := \sum_{n \leq x} (f \ast 1)(n) = \sum_{d \leq x} f(d) \Floor{x}{d}, 
\]
inherit partition-function-like recurrence relations from the structure of the factorizations we have 
constructed. In particular, for $n, x \geq 1$ we have that \cite{AA} 
\begin{align*} 
(f \ast 1)(n+1) & = \sum_{b = \pm 1} \sum_{k=1}^{\left\lfloor \frac{\sqrt{24n+1}-b}{6}\right\rfloor} 
(-1)^{k+1} (f \ast 1)\left(n+1-\frac{k(3k+b)}{2}\right) + 
\sum_{k=1}^{n+1} s_{n+1,k} f(k), \\ 
\Sigma_f(x+1) & = \sum_{b = \pm 1} \sum_{k=1}^{\left\lfloor \frac{\sqrt{24x+1}-b}{6}\right\rfloor} 
(-1)^{k+1} \Sigma_f\left(n+1-\frac{k(3k+b)}{2}\right) + 
\sum_{n=0}^x \sum_{k=1}^{n+1} s_{n+1,k} f(k).
\end{align*} 

\subsubsection{Generalized Lambert series expansions} 

Most generally in \cite{MERCA-SCHMIDT-RAMJ}, we define the generalized Lambert series factorizations 
which are parameterized by the lower triangular 
sequence on the right-hand-side of the next equation by 
the following series expansion 
for integers $0 \leq \beta < \alpha$, 
$C(q)$ any convergent OGF for $|q| < 1$ such that $C(0) \neq 0$, and $\bar{f}$ a function depending on $f$: 
\begin{equation} 
\label{eqn_GenFactThmExp_def_v1} 
L_f(\alpha, \beta; q) := \sum_{n \geq 1} \frac{f(n) q^{\alpha n-\beta}}{1-q^{\alpha n-\beta}} = 
     \frac{1}{C(q)} \sum_{n \geq 1} \left(\sum_{k=1}^{n} \bar{s}_{n,k}(\alpha, \beta) \bar{f}(k)\right) q^{n}; |q| < 1. 
\end{equation}
The lower triangular sequence, $\bar{s}_{n,k}(\alpha, \beta)$, 
in the expansion above is invertible and depends on the parameters  
$(\alpha, \beta, f, C(q))$ 
that define this Lambert series expansion. 
Alternately, for $|q|<1$ and integers $0\leq \beta<\alpha$, 
we have that 
	\begin{align*}
	\sum_{n=1}^{\infty} a_n \frac{q^{\alpha n-\beta}}{1-q^{\alpha n-\beta}} = 
	\frac{1}{(q^{\alpha-\beta};q^\alpha)_\infty} \sum_{n=1}^{\infty} \left( \sum_{k=1}^n (s_o(n,k)-s_e(n,k)) a_k\right)  q^n,
	\end{align*}
	where $s_o(n,k)$ and $s_e(n,k)$ denote the number of $(\alpha k-\beta)$'s in all partitions of $n$ into an odd 
	(respectively even) number of distinct parts of the form $\alpha k-\beta$.
Similarly, for $|q|<1$, $0\leq \beta<\alpha$, we can expand  
	\begin{align*}
    %\label{eqn_FactThmExp_def_v3}
	\sum_{n=1}^{\infty} a_n \frac{q^{\alpha n-\beta}}{1-q^{\alpha n-\beta}} = 
	(q^{\alpha-\beta};q^\alpha)_\infty \sum_{n=1}^{\infty} \left( \sum_{k=1}^n s(n,k) a_k\right)  q^n,
	\end{align*}
	where $s(n,k)$ denotes the number of $(\alpha k-\beta)$'s in all partitions of $n$ into parts of the form $\alpha k-\beta$.
Moreover, for any fixed arithmetic function $f$, 
if we define the arithmetic functions $\gamma(n)$ and 
$\widetilde{\gamma}(n) := \sum_{d|n} \gamma(d)$ that depend on 
a fixed factorization pair $(C(q), \bar{s}_{n,k})$ in \eqref{eqn_GenFactThmExp_def_v1} such that 
\begin{align*}
\bar{s}_{n,k}^{(-1)} & := \sum_{d|n} [q^{d-k}] \frac{1}{C(q)} \cdot \gamma(n/d), 
\end{align*} 
then we have that the sequence of $\bar{f}(n)$ is given by the following formula for all $n \geq 1$: 
\begin{align*} 
\bar{f}(n) & = 
     \sum_{\substack{d|n \\ d \equiv \beta \bmod \alpha}} 
     f\left(\frac{d-\beta}{\alpha}\right) 
     \widetilde{\gamma}\left(\frac{n}{d}\right). 
\end{align*} 

\newpage
\section{Ordinary Lambert series (LGF) identities} 

\subsection{Listings of identities for  
            arithmetic functions} 

We have the following well-known ``classical'' examples of Lambert series identities 
\cite[\S 27.7]{NISTHB} \cite[\S 17.10]{HARDYWRIGHT} \cite[\S 11]{APOSTOLANUMT}:
\begin{subequations}
\begin{align} 
\label{eqn_WellKnown_LamberSeries_Examples} 
\sum_{n \geq 1} \frac{\mu(n) q^n}{1-q^n} & = q, \\ 
\label{eqn_WellKnown_LamberSeries_Examples_phi} 
\sum_{n \geq 1} \frac{\phi(n) q^n}{1-q^n} & = \frac{q}{(1-q)^2}, \\ 
\sum_{n \geq 1} \frac{n^{\alpha} q^n}{1-q^n} & =  
     \sum_{m \geq 1} \sigma_{\alpha}(n) q^n, \\ 
\sum_{n \geq 1} \frac{\lambda(n) q^n}{1-q^n} & = \sum_{m \geq 1} q^{m^2}, \\ 
\sum_{n \geq 1} \frac{\Lambda(n) q^n}{1-q^n} & = \sum_{m \geq 1} \log(m) q^m, \\ 
\sum_{n \geq 1} \frac{|\mu(n)| q^n}{1-q^n} & = \sum_{m \geq 1} 2^{\omega(m)} q^m, \\ 
\sum_{n \geq 1} \frac{J_t(n) q^n}{1-q^n} & = \sum_{m \geq 1} m^t q^m, \\ 
\sum_{n \geq 1} \frac{\mu(\alpha n) q^n}{1-q^n} & = - \sum_{n \geq 0} q^{\alpha^n}, \alpha \in \mathbb{P} \\ 
\sum_{n \geq 1} \frac{q^n}{1-q^n} & = \frac{\psi_q(1) + \log(1-q)}{\log(q)}, \\ 
\sum_{n \geq 1} \frac{\operatorname{lsb}(n) q^n}{1-q^n} & = \frac{\psi_{q^2}(1/2) + \log(1-q^2)}{2 \log(q)}.
\end{align}
\end{subequations}

\subsection{Other LGF-variant identities} 

For any arithmetic function $f$ and integers $k \geq 1$, we have that 
\begin{align} 
\sum_{n \geq 1} \frac{f(n) q^{n^k}}{1-q^{n^k}} & = \sum_{m \geq 1} \left(\sum_{d^k | n} f(d)\right) q^m.
\end{align} 
This follows from the uniform expansion in the last equation as 
\begin{align*}
\sum_{n \geq 1} \sum_{m \geq 1} f(n) q^{m n^k}. 
\end{align*}
The integers $r$ in the series expansion of the 
last equation such that we have $f(r)$ as a coefficient of 
$q^p$ correspond to the divisors $r$ of $p = m \cdot n^k$. 

\begin{subequations} 
For example, we have that 
\begin{align}
\sum_{n \geq 1} \frac{\lambda_k(n) q^n}{1-q^n} & = \sum_{m \geq 1} q^{m^k} \\ 
\sum_{n \geq 1} |\mu_k(n)| q^n & = q^{k+1} \\ 
\sum_{n \geq 1} \frac{\mu(n) q^{n^2}}{1-q^{n^2}} & = \sum_{m \geq 1} |\mu(m)| q^m. 
\end{align} 
\end{subequations} 

\section{Modified Lambert series identities} 

\subsection{Definitions} 

For $|q| < 1$ and $f$ any arithmetic function, let 
\[
\widehat{L}_f(q) := \sum_{n \geq 1} \frac{f(n) q^n}{1+q^n}. 
\]
We can write $\widehat{L}_f(q) \equiv L_h(q)$, e.g., as an ordinary Lambert series expansion, where 
\[
h(n) = \begin{cases} 
     h(n), & \text{ if $n$ is odd; } \\ 
     h(n) - 2h\left(\frac{n}{2}\right), & \text{ if $n$ is even. } 
     \end{cases}
\]
Thus we have that 
\begin{equation} 
\label{eqn_} 
\widehat{L}_f(q) = L_f(q) - 2L_f(q^2). 
\end{equation} 

\subsection{Examples} 

The two primary Lambert series expansions from the previous section that admit ``nice'', algebraic 
closed-form expressions are translated below:
\begin{subequations} 
\begin{align} 
\sum_{n \geq 1} \frac{\mu(n) q^n}{1+q^n} & = q-2q^2 \\ 
\sum_{n \geq 1} \frac{\phi(n) q^n}{1+q^n} & = \frac{q(q+q^2)}{(1-q^2)^2}. 
\end{align} 
\end{subequations} 
Note that Section \ref{subSection_GCDSums} also provides a pair of related series in the context of 
GCD sums of an arithmetic function. 

\section{Generalized Lambert series identities} 

\subsection*{Listings of identities} 

From \cite[\S 17.10]{HARDYWRIGHT}, we obtain that 
\begin{equation}
\sum_{n \geq 1} \frac{4 \cdot (-1)^{n+1} q^{2n+1}}{1-q^{2n+1}} = \sum_{m \geq 1} r_2(m) q^m, |q| < 1. 
\end{equation} 
\begin{definition}[Jacobi theta functions]
For complex-valued $q, z$, we define the next four 
variants of the \emph{Jacobi theta functions} by the 
following bilateral series:
\begin{align*} 
\vartheta_1(z, q) & := \sum_{n=-\infty}^{\infty} 
     (-1)^{n+1/2} q^{(n+1/2)^2} e^{(2n+1)\imath z}, \\ 
\vartheta_2(z, q) & := \sum_{n=-\infty}^{\infty} 
     q^{(n+1/2)^2} e^{(2n+1)\imath z}, \\ 
\vartheta_3(z, q) & := \sum_{n=-\infty}^{\infty} 
     q^{n^2} e^{2n \imath z}, \\ 
\vartheta_4(z, q) & := \sum_{n=-\infty}^{\infty} 
     (-1)^{n} q^{n^2} e^{2n \imath z}
\end{align*}
\end{definition}

There are a number of classical theta function related series of the following forms \cite[\S 20]{NISTHB}: 
\begin{subequations} 
\begin{align} 
\sum_{n \geq 1} \frac{q^n}{1+q^{2n}} & = \frac{1}{4}\left[\vartheta_3^2(q) - 1\right] \\ 
\sum_{n \geq 1} \frac{q^{2n+1}}{1+q^{4n+2}} & = \frac{1}{4}\left[\vartheta_3^2(q) - \vartheta_2^2(q)\right] = 
     \frac{\vartheta_2^2(q^2)}{4} \\ 
\sum_{n \geq 1} \frac{q^n}{1-q^{2n}} & = L_1(q) - L_1(q^2) \\ 
\sum_{n \geq 1} \frac{q^{2n+1}}{1-q^{4n+2}} & = L_1(q) - 2L_1(q^2) + L_1(q^4) \\ 
\sum_{n \geq 1} \frac{4 \sin(2nz) q^{2n}}{1-q^{2n}} & = \frac{\vartheta_1^{\prime}(z, q)}{\vartheta_1(z, q)} \\ 
\sum_{n \geq 1} \frac{4 (-1)^n \sin(2nz) q^{2n}}{1-q^{2n}} & = \frac{\vartheta_2^{\prime}(z, q)}{\vartheta_2(z, q)} \\ 
\sum_{n \geq 1} \frac{4 (-1)^n \sin(2nz) q^{n}}{1-q^{2n}} & = \frac{\vartheta_3^{\prime}(z, q)}{\vartheta_3(z, q)} \\ 
\sum_{n \geq 1} \frac{4 \sin(2nz) q^{n}}{1-q^{2n}} & = \frac{\vartheta_4^{\prime}(z, q)}{\vartheta_4(z, q)} \\ 
\sum_{n=-\infty}^{\infty} \frac{(-1)^n e^{2\imath nz} q^{n^2}}{q^{-n} e^{-\imath z} + q^n e^{\imath z}} & = 
     \frac{\vartheta_2(0, q) \vartheta_3(z, q) \vartheta_4(z, q)}{\vartheta_2(z, q)}. 
\end{align}
\end{subequations} 
There are a number of Lambert, and Lambert-like series, for mock theta functions of order $6$ given in 
\cite[\S 8]{RAMANUJAN-LSERIES-SURVEY}. These series expansions are cited as follows where 
$J_{a,m} := (q^a, q^{m-a}, q^m; q^m)_{\infty}$: 
\begin{subequations} 
\label{eqn_MockThetaSixthOrder_idents} 
\begin{align} 
\phi_{\mock}(q) & = \sum_{n \geq 0} \frac{(-1)^n q^{n^2} (q; q^2)_n}{(-q)_{2n}} && = 
     \frac{2}{J_{1,3}} \sum_{r=-\infty}^{\infty} \frac{q^{r(3r+1)/2}}{1+q^{3r}} \\ 
\Psi_{\mock}(q) & = \sum_{n \geq 0} \frac{(-1)^n q^{(n+1)^2} (q; q^2)_n}{(-q)_{2n+1}} && = 
     \frac{2}{J_{1,3}} \sum_{r=-\infty}^{\infty} \frac{q^{r(3r+1)/2}}{1+q^{3r+1}} \\ 
\rho_{\mock}(q) & = \sum_{n \geq 0} \frac{q^{n(n+1)/2} (-q)_n}{(q; q^2)_{n+1}} && = 
     \frac{1}{J_{1,6}} \sum_{r=-\infty}^{\infty} \frac{(-1)^r q^{r(3r+4)}}{1-q^{6r+1}} \\ 
\sigma_{\mock}(q) & = \sum_{n \geq 0} \frac{q^{n(n+2)/2} (-q)_n}{(q; q^2)_{n+1}} && = 
     \frac{1}{J_{1,6}} \sum_{r=-\infty}^{\infty} \frac{(-1)^r q^{(r+1)(3r+1)}}{1-q^{6r+3}} \\ 
\gamma_{\mock}(q) & = \sum_{n \geq 0} \frac{q^{n^2} (q)_n}{(q^3; q^3)_{n}} && = 
     \frac{1}{(q; q)_{\infty}} \sum_{r=-\infty}^{\infty} \frac{(-1)^r q^{r(3r+1)/2}}{1+q^r+q^{2r}}. 
\end{align} 
\end{subequations} 

\section{Lambert series over Dirichlet convolutions and Apostol divisor sums} 

\subsection{Dirchlet convolutions} 

\subsubsection{Definitions} 

Given two prescribed arithmetic functions $f$ and $g$ we define their \emph{Dirichlet convolution}, 
denoted by $h = f \ast g$, to be the function 
\begin{align*} 
(f \ast g)(n) & := \sum_{d|n} f(d) g\left(\frac{n}{d}\right), 
\end{align*} 
for all natural numbers $n \geq 1$ \cite[\S 2.6]{APOSTOLANUMT}. The classical M\"obius inversion 
result is stated in terms of convolutions as follows, where 
$\mu$ is the M\"obius function: $h = f \ast 1$ if and only if $f = h \ast \mu$. 
There is a natural connection between the coefficients of the Lambert series of an arithmetic function $a_n$ and its corresponding \emph{Dirichlet generating function}, 
$\DGF(a_n; s) := \sum_{n \geq 1} a_n / n^s$. 
Namely, we have that for any $s \in \mathbb{C}$ such that $\Re(s) > 1$ 
$$b_n = [q^n] \sum_{n \geq 1} \frac{a_n q^n}{1-q^n} \quad\text{ if and only if }\quad 
  \DGF(b_n; s) = \DGF(a_n; s) \zeta(s), $$ where $\zeta(s)$ is the Riemann zeta function. Moreover, we can further connect the coefficients of the Lambert series over a convolution of arithmetic functions to its 
associated Dirichlet series by noting that $\DGF(f \ast g; s) = \DGF(f; s) \cdot \DGF(g; s)$. 

\begin{notation}[Expanding Dirichlet inverse functions]
The \emph{Dirichlet inverse function} of $f(n)$, denoted $f^{-1}(n)$, is an arithmetic function such that 
$(f \ast f^{-1})(n) = \delta_{n,1}$ for all $n \geq 1$. The function $f^{-1}$ exists and is unique if and only if 
$f(1) \neq 0$. In these cases, we can expand the inverse function in terms of weighted terms in $f$ recursively 
according to the formula 
\[
f^{-1}(n) = \begin{cases} 
     \frac{1}{f(1)}, & n = 1; \\ 
     -\frac{1}{f(1)} \times \sum\limits_{\substack{d|n \\ d>1}} f(d) f^{-1}\left(\frac{n}{d}\right), & n \geq 2. 
     \end{cases} 
\]
We have that \cite{MOUSAVI-SCHMIDT-2019} 
\[
f^{-1}(n) = \sum_{j=1}^{\Omega(n)} \frac{(-1)^j \cdot (f-f(1)\varepsilon)_{\ast_j}(n)}{f(1)^{j+1}}. 
\]
Note that Section \ref{Section_OtherSpecialIdents} contains a formula enumerating the Dirichlet inverse of 
any Dirichlet invertible arithmetic function $f$. 
\end{notation}

\subsubsection{General identities} 

We can see that the Lambert series over the convolution $(f \ast g)(n)$ is given by the double sum 
\[
L_{f \ast g}(q) = \sum_{n \geq 1} f(n) L_g(q^n), |q| < 1. 
\]
Similarly, 
\[
\widehat{L}_{f \ast g}(q) = \sum_{n \geq 1} f(n) \left[L_g(q^n)-2L_g(q^{2n})\right]. 
\]
Clearly we have by M\"obius inversion that the \emph{ordinary generating function} (OGF) of $f$ is 
given by 
\[
L_{f \ast \mu}(q) = \sum_{n \geq 1} f(n) q^n. 
\]
If $F(x) := \sum_{n \leq x} f(n)$ is the \emph{summatory function} of $f$, then we have that 
\[
\sum_{n \geq 1} F(n) q^n = \sum_{n \geq 1} \mu(n) \frac{L_f(q^n)}{1-q}. 
\]
\begin{proof}
The last identity follows by writing 
\begin{align*}
[q^n] \frac{L_f(q)}{1-q} & = \sum_{k \leq n} (f \ast 1)(k) \\ 
     & = \sum_{k=1}^n F(k) \sum_{r=\Floor{n}{k+1}+1}^{\Floor{n}{k}} 1 \\ 
     \implies [q^n] L_f(q) & = \sum_{m|n} \frac{n}{m} \times 
     (F(m)-F(m-1)). 
\end{align*} 
Thus it follows that since $\operatorname{Id}_k^{-1}(n) = \mu(n) \operatorname{Id}_k(n) = \mu(n) n^{k}$
as in \cite[\cf \S 2]{APOSTOLANUMT}, we get that 
\[
\frac{L_{f \ast \operatorname{Id}_1^{-1}}(q)}{1-q} = \sum_{n \geq 1} \mu(n) \frac{L_f(q^n)}{1-q} = 
     \sum_{n \geq 1} F(n) q^n. 
     \qedhere 
\]
\end{proof} 

\subsubsection{Listing of particular identities} 

We have the following convolution identities for Lambert series expansions of special functions 
\cite[\S 7.4]{GKP} \cite[\S 24.4(iii)]{NISTHB}: 
\begin{subequations}
\begin{align} 
\sum_{n \geq 1} \frac{\psi_k(n) q^n}{1-q^n} & = \sum_{j=0}^{k} \gkpSII{k}{j} j! \times 
     \sum_{m \geq 1} 2^{\omega(m)} \frac{q^{mj}}{(1-q^m)^{j+1}} \\ 
\notag 
     & = \sum_{j=0}^{k} \gkpSII{k}{j} j! \times 
     \sum_{n \geq 1} \sum_{d|\lfloor n/j\rfloor} 2^{\omega(d)} \binom{\Floor{n}{j} \frac{1}{d} + j}{j} \cdot q^n; 
     m \in \mathbb{N}, \\ 
\sum_{n \geq 1} \frac{(\sigma_k \ast \mu)(n) q^n}{1-q^n} & = \sum_{j=0}^{k} \gkpSII{k}{j} \frac{j! q^j}{(1-q)^{j+1}}; 
     m \in \mathbb{N}; \sigma_k \ast \mu = \operatorname{Id}_k, \\ 
\sum_{n \geq 1} \frac{\sigma_1(n) q^n}{1-q^n} & = \sum_{n \geq 1} \frac{d(n) q^n}{(1-q^n)^2}; \sigma_1 = \phi \ast \sigma_0, \\ 
\sum_{n \geq 1} \frac{(\phi_k \ast \operatorname{Id}_k)(n) q^n}{1-q^n} & = 
     \sum_{m \geq 1} \frac{1}{k+1} \times 
     \left(B_{k+1}(m+1) - B_{k+1}(0)\right) q^m. 
\end{align}
\end{subequations} 
\begin{subequations}
We also have some unique generating function expressions for common summatory functions, including the 
following identity:
\begin{equation}
\sum_{n \geq 1} \mu(n) \frac{L_{\mu \ast \omega}(q^n)}{1-q} = \sum_{x \geq 1} \pi(x) q^x. 
\end{equation}
\end{subequations}

\subsubsection{Characteristic functions} 

\begin{subequations}
\label{eqn_DirCvlLSeries_CharFunc_Idents} 
The argument used to arrive at the last identity shows that if $A \subseteq \mathbb{Z}^{+}$ and its 
indicator function is denoted by $\chi_A(n)$, then we have that 
\begin{equation}
\sum_{n \geq 1} \mu(n) L_{\chi_A}(q^n) = \sum_{a \in A} q^a. 
\end{equation}
For example, if $\mathbb{N}_{\operatorname{sqfree}}$ denotes the set of positive squarefree 
integers, then 
\begin{equation}
\sum_{n \geq 1} \mu(n) L_{\mu^2}(q^n) = \sum_{k \in \mathbb{N}_{\operatorname{sqfree}}} q^k. 
\end{equation} 
Moroever, if $\chi_A(n) = (\mu \ast g_A)(n)$, then 
\begin{equation} 
\sum_{n \geq 1} \frac{g_A(n) f(n) q^n}{1-q^n} = \sum_{a \in A} L_f(q^a). 
\end{equation} 
For example, in equation \eqref{eqn_PointwiseProductIdentsOmegaFuncs} 
of the next section, we prove a 
prime summation identity for the Lambert series over the pointwise products of $\omega(n) f(n)$ 
and $\lambda(n) f(n)$ for any arithmetic $f$. 
\end{subequations} 

\subsubsection{Expressions for series generating Dirichlet inverse functions} 

We denote by $f_{\ast_j}(n)$ the $j$-fold convolution of $f$ with itself, i.e., the sequence 
defined recursively by 
\[
f_{\ast_j}(n) = \begin{cases} 
     \delta_{n,1}, & \text{ if $j=0$; } \\ 
     \sum\limits_{d|n} f(d) f_{\ast_{(j-1)}}\left(\frac{n}{d}\right), & 
     \text{ if $j \geq 1$. } 
     \end{cases}
\]
Then as in \cite{MOUSAVI-SCHMIDT-2019}, we have that 
\[
f^{-1}(n) = \sum_{j=1}^{\Omega(n)} \binom{\Omega(n)}{j} \frac{(-1)^j}{f(1)^{j+1}} f_{\ast_j}(n). 
\]
Hence, we have that 
\begin{align} 
\sum_{n \geq 1} \frac{f^{-1}(n) q^n}{1-q^n} = \sum_{n \geq 1} \left(1-\frac{f(n)}{f(1)}\right)^{\Omega(n)} \times 
     \frac{L_f(q^n)}{f(1)f(n)} - \sum_{n \geq 1} \frac{L_f(q^n)}{f(1)f(n)}. 
\end{align}

\subsection{GCD transform sums} 
\label{subSection_GCDSums} 

\subsubsection{General identities} 

We have that 
\begin{subequations}
\begin{align} 
\label{eqn_GCDTransGenIdents_v1} 
\sum_{n \geq 1} \left(\sum_{\substack{1 \leq d \leq n \\ (d, n)=1}} f(d)\right) \frac{q^n}{1-q^n} & = 
     \sum_{k \geq 1} \left(\sum_{d|k} \frac{\mu(d)}{1-q^d}\right) f(k) q^k \\ 
\label{eqn_GCDTransGenIdents_v2} 
\sum_{n \geq 1} \left(\sum_{\substack{1 \leq d \leq n \\ (d, k)=m}} f(d)\right) \frac{q^n}{1-q^n} & = 
     \sum_{k \geq 1} \left(\sum_{d|k} \frac{\mu(d)}{1-q^{md}}\right) f(k) q^k \\ 
\label{eqn_GCDTransGenIdents_v3} 
\sum_{n \geq 1} \left(\sum_{d=1}^{n} f(\gcd(d, n))\right) \frac{q^n}{1-q^n} & = 
     \sum_{n \geq 1} \frac{(f \ast \phi)(n) q^n}{1-q^n} \\ 
\notag 
     & = \sum_{n \geq 1} f(n) \frac{q^n}{(1-q^n)^2} \\ 
\notag
     & = \sum_{n \geq 1} \sum_{k=1}^n (f \ast 1)(\gcd(k, n)) q^n. 
\end{align}
\end{subequations}

The identities in \eqref{eqn_GCDTransGenIdents_v1} and 
\eqref{eqn_GCDTransGenIdents_v2} result from the 
following equation for fixed integers $k \geq 1$ and 
$1 \leq m \leq k$ \citep[\S 3.2]{MOUSAVI-SCHMIDT-2019}: 
\[
\sum_{n \geq k} \Iverson{(n, k) = m} q^n = 
     \sum_{d|k} \frac{q^k \mu(d)}{1 - q^{md}}. 
\]
The results in \eqref{eqn_GCDTransGenIdents_v3} are respectively 
consequences of the following equation, 
\eqref{eqn_WellKnown_LamberSeries_Examples_phi} and the 
Dirichlet convolution identity that $(\phi \ast 1)(n) = n$ for 
all $n \geq 1$
\citep{KAMP-GCD-TRANSFORM} \citep[\S 27.5]{NISTHB}: 
\[
\sum_{1 \leq k \leq m} h(\gcd(k, m)) = (h \ast \phi)(m), 
     m \geq 1. 
\]

\subsubsection{Particular cases} 

In \cite{KAMP-GCD-TRANSFORM}, formulas for the discrete Fourier transform of a function 
evaluated at a $\gcd$ argument are derived. The reference 
also connects Lambert series expansions of Liouville 
for the divisor sum functions $\phi_a(n)$ 
(non-standard notation) that generalize the classical \emph{Euler totient function} as 
\begin{subequations} 
\begin{align} 
\sum_{n \geq 1} \left(\sum_{d|(a,n)} d\cdot 
     \phi\left(\frac{n}{d}\right)\right) \frac{q^n}{1-q^n} & = 
     \frac{\sum\limits_{k=1}^{2a} (a-|k-a|) d(\gcd(a-|k-a|, a)) q^k}{(1-q^a)^2} \\ 
\sum_{n \geq 1} \left(\sum_{d|(a,n)} d\cdot 
     \phi\left(\frac{n}{d}\right)\right) \frac{q^n}{1+q^n} & = 
     \frac{p[a](q)}{(1-q^{2a})^2}, 
\end{align} 
where 
\[
p[a](q) := \sum\limits_{k=1}^{4a} 
     \left[(2a-|k-2a|) d(\gcd(2a-|k-2a|, a)) - \Iverson{k \text{\ even}} \left(a-|k/2-a|\right) 
     d(\gcd(a-|k/2-a|, a))\right] q^k. 
\]
\end{subequations} 

There are related LCM Dirichlet series, or DGF, identities that we can cite to find a 
Lambert series expansion for these functions from \cite{CATALOG-INTDIRSERIES}: 
\begin{subequations} 
\begin{align} 
\sum_{n \geq 1} \left(\sum_{k=1}^{n} [k,n]\right) \frac{q^n}{1-q^n} & = 
     \sum_{m \geq 1} \frac{1}{2} \left(\sigma_1(m) + \sum_{d|m} \sum_{r|\frac{m}{d}} 
     d \sigma_2(d) \mu\left(\frac{m}{dr}\right) \left(\frac{m}{dr}\right)^2\right) q^m, \\ 
\sum_{n \geq 1} \left(\sum_{k=1}^{n} [k,n]^m\right) \frac{q^n}{1-q^n} & = \sum_{n \geq 1} \left(
     \sigma_m(n) + \sum_{i=1}^{m+1} \binom{m+1}{i} \frac{B_{m+1-i}}{m+1} 
     \left(1 \ast \operatorname{Id}_m \ast \operatorname{Id}_{m+i} \ast \operatorname{Id}_{2m}^{-1}\right)(n)\right) 
     q^m. 
\end{align} 
\end{subequations} 
In particular, we learn from \cite[D-19; D-71]{CATALOG-INTDIRSERIES} that the first Lambert series 
above is generated as $L_{f_1}(q)$ when 
$f_1 = \frac{1}{2}\left(\operatorname{Id}_1 \ast \operatorname{Id}_2^{-1} \cdot (\operatorname{Id}_2+\operatorname{Id}_3)\right)$, 
and the second (LCM powers sum) series is generated as $L_{f_2}(q)$ with 
$$f_2 = \operatorname{Id}_m + \sum_{i=1}^{m+1} \binom{m+1}{i} \frac{B_{m+1-i}}{m+1} \left( 
 \operatorname{Id}_m \ast \operatorname{Id}_{m+i} \ast \operatorname{Id}_{2m}^{-1}\right).$$ 

\subsection{Anderson-Apostol divisor sums} 

We consider the Lambert series generating functions 
over the sums \citep{MOUSAVI-SCHMIDT-2019}
\[
S_{1,m}(f, g; n) := \sum_{d|(m,n)} f(d) g\left(\frac{m}{d}\right). 
\]
We have that 
\begin{align} 
\notag 
\widetilde{L}_{1,m}(f, g; q) & := 
     \sum_{n \geq 1} \frac{S_{1,m}(f, g; n) q^n}{1-q^n} \\ 
\label{eqn_ApostolSumsIdent_PrimaryFormula_v1} 
     & \phantom{:}= \sum_{n \geq 1} (f \ast g \ast 1)(\gcd(m, n)) q^n. 
\end{align} 

\begin{proof}[Proof of \eqref{eqn_ApostolSumsIdent_PrimaryFormula_v1}] 
We prove the following for integers $n \geq 1$:
\begin{align*} 
\sum_{d|n} \sum_{r|(m,d)} f(r) g\left(\frac{m}{r}\right) & = 
     \sum_{\substack{s|m \\ s|d|n}} \sum_{r|s} 
     f(r) g\left(\frac{s}{r}\right) \\ 
     & = \sum_{s|(m,n)} (f \ast g)(s). 
\end{align*} 
The key transition step in the above equations is in noting that $(d, m)$ is a divisor of 
both $d$ and $m$ for any integers $d,m \geq 1$. 
\end{proof} 

The primary special case of interest with these types of sums is the \emph{Ramanujan sum}, 
$c_q(x) \equiv S_{1,x}(\operatorname{Id}_1, \mu; q)$. 
In particular, as expanded in \citep{KAMP-GCD-TRANSFORM}, we know that 
\[
c_q(x) = \sum_{d|(q, x)} d \mu\left(\frac{k}{d}\right), 
     \text{ for\ integers\ } q, x \geq 1. 
\]
These sums are periodic modulo 
$m \geq 1$ and have a finite Fourier series expansion with known coefficients. 
Let 
\[
a_k(f, g; m) = \sum_{d|(m,k)} g(d) f\left(\frac{k}{d}\right) \cdot 
     \frac{d}{k}.
\]
Then we have that \cite[\S 27.10]{NISTHB} 
\[
S_{1,m}(f, g; n) = \sum_{k=1}^m a_m(f, g; m) \cdot e^{2\pi\imath \cdot kn/m}. 
\]

\subsection{Another summation variant} 

We next consider the Lambert series generating functions
over the sums 
\[
S_{2,m}(f, g; n) := \sum_{d|(m,n)} f(d) g\left(\frac{mn}{d^2}\right). 
\]
As an example of an identity involving this summation type, 
we have that the \emph{Ramanujan tau function}, $\tau(n)$, 
satisfies 
\[
\tau(m) \tau(n) = \sum_{d|(m, n)} d^{11} \tau\left(\frac{mn}{d^2}\right). 
\]
Additionally, as another example, 
for any $\alpha \in \mathbb{C}$ and $m, n \geq 1$, 
\[
\sigma_{\alpha}(m) \sigma_{\alpha}(n) = \sum_{d|(m,n)} d^{\alpha} \sigma_{\alpha}\left(\frac{mn}{d^2}\right). 
\]
Note that if $g$ is completely multiplicative, then \cite[\S 2, Exercise 31]{APOSTOLANUMT} 
\[
f(m) f(n) = \sum_{d|(m, n)} g(d) f\left(\frac{mn}{d^2}\right). 
\]

\section{Other special identities} 
\label{Section_OtherSpecialIdents} 

\subsection{Hadamard products with special arithmetic functions} 

We have by Mobius inversion and our previous identities on Dirichlet convolutions that: 
\begin{subequations} 
\label{eqn_PointwiseProductIdentsOmegaFuncs} 
\begin{align} 
\label{eqn_HPOtherSpecIdents_v1} 
\sum_{n \geq 1} \frac{\omega(n) f(n) q^n}{1-q^n} & = \sum_{p\text{\ prime}} L_f(q^p) \\ 
\label{eqn_HPOtherSpecIdents_v2} 
\sum_{n \geq 1} \frac{\lambda(n) f(n) q^n}{1-q^n} & = \sum_{d \geq 1} \sum_{n \geq 1} \frac{\mu(n) f(nd^2) q^{nd^2}}{1-q^{nd^2}}. 
\end{align} 
\end{subequations}

\begin{proof} 
These two equations follow from \eqref{eqn_DirCvlLSeries_CharFunc_Idents} by noting that the 
characteristic function of the primes is given by $\chi_{\mathbb{P}} = \omega \ast \mu$ and that 
the characteristic function of the squares is given by 
$\chi_{\operatorname{sq}} = \lambda \ast \mu$. 
\end{proof} 

\subsection{Results on divisor sums involving products of $\omega(n)$ and $\mu(n)$} 

\begin{subequations}
Suppose that $f$ is multiplicative such that $f(p) \neq +1, -1$, respectively, for all primes $p$. Then we have that 
\cite{TWAKHARE-2016} 
\begin{align} 
\sum_{d|n} \mu(d) \omega(d) f(d) & = \prod_{p|n} (1-f(p)) \times \sum_{p|n} \frac{f(p)}{f(p)-1} \\ 
\sum_{d|n} |\mu(d)| \omega(d) f(d) & = \prod_{p|n} (1+f(p)) \times \sum_{p|n} \frac{f(p)}{1+f(p)}
\end{align} 
\end{subequations} 
\begin{subequations} 
Under the same respective conditions, suppose that $f$ is indeed completely multiplicative. 
Then similarly, we obtain that 
\begin{align} 
\sum_{d|n} \mu(d) \omega(d) f(d) & = \sum_{d|n} \mu(d) f(d) \times \sum_{p|n} \frac{f(p)}{f(p)-1} \\ 
\sum_{d|n} |\mu(d)| \omega(d) f(d) & = \sum_{d|n} |\mu(d)| f(d) \times \sum_{p|n} \frac{f(p)}{1+f(p)}
\end{align} 
\end{subequations} 

\subsection{Divisor sum convolution identities involving other prime-related arithmetic functions} 

We have the following prime sum related divisor sum identities in the form of $f \ast 1$ generated 
by a Lambert series generating function over a multiplicative $f$ \cite{BC-IDENTS-MNUMT}: 
\begin{subequations} 
\begin{align} 
\sum_{d|n} \frac{\mu(d) \log d}{d} & = \frac{\phi(n)}{n} \sum_{p|n} \frac{\log p}{1-p}, \\ 
\sum_{d|n} \frac{|\mu(d)| \log d}{d^k} & = \frac{\Psi_k(n)}{n^k} \sum_{p | n} \frac{\log p}{p^{k}+1}, \\ 
\sum_{d|n} \frac{|\mu(d)| \log d}{\phi(d)} & = \frac{n}{\phi(n)} \sum_{p|n} \frac{\log p}{p}, \\ 
\sum_{d|n} \frac{\mu(d) \log d}{\sigma_0(d)} & = -2^{\omega(n)} \log \gamma(n), \\ 
\sum_{d|n} \mu(d)^{a} d_k(d) \log d & = (1+(-1)^a k)^{\omega(n)} \times 
     \frac{k \log \gamma(n)}{k + (-1)^a}; a \in \{1,2\}, k \geq 2, \\ 
\sum_{d|n} \mu(d) \sigma_1(d) \log d & = (-1)^{\omega(n)} \gamma(n) \left( 
     \log \gamma(n) + \sum_{p|n} \frac{\log p}{p}\right), \\ 
\sum_{d|n} \mu(d)^{a} f(d) \log d & = \prod_{p|n} (1+(-1)^{a} f(p)) \times 
     \sum_{p|n} \frac{f(p) \log p}{f(p)+(-1)^{a}}; a \in \{1,2\}, \\ 
\sum_{d|n} |\mu(d)| k^{\omega(d)} & = (k+1)^{\omega(n)}. 
\end{align} 
\end{subequations}

\subsection{Relations of generalized Lambert series to $q$-series expansions} 
\label{subSection_OtherProps_RefsToq-SeriesExps} 

We do not focus on connections of other forms of generalized Lambert series expansions 
to $q$-series and partition generating functions, nor consider their representations in the 
context of modular forms. 
In this sense, we note that one can consider a class of generalized Lambert series defined by 
\[
L(\alpha; t, q) := \sum_{n \geq 1} \frac{t^n}{1-x q^n},  
\]
and then connect variants of this function to $q$-series (see, for example, the identities given in 
\eqref{eqn_MockThetaSixthOrder_idents}). 
For an overview of that vast material, 
we refer the reader to a subset of relevant references in 
\cite{ARNDT-GENLAMBERTSERIES,ANDREWS-BAILEYCH-LAMSERIES-1992,RAMANUJAN-LSERIES-SURVEY}. 

\nocite{SCHMIDT-GTDISSERTATION}

\newpage
\renewcommand{\refname}{References} 
%\bibliography{thesis-references}{}
\bibliographystyle{plain}

\end{document}